\documentclass[a4paper,12pt,twoside]{article}

   \usepackage  [pics]{optional} 

\makeatletter

\renewcommand\section{\@startsection 
  {section}{1}{\z@}%
  {-3.25ex\@plus -1ex \@minus -.2ex}
  {1.5ex \@plus .2ex}%
  {\normalfont\large\myem}}
                                                                    
\makeatother

\usepackage{amsthm}
\usepackage{amsmath}
\usepackage{amssymb}
\usepackage{bbm}
\usepackage{pstricks,pst-node,multido}
\usepackage{ragged2e}

\renewcommand{\DeclareMathOperator}[2]
{\newcommand{#1}{\mathop{\ensuremath{\text{\upshape #2}}}}}

\newcommand{\normalset}[2] {\{#1\ |\ #2\}}
\newcommand{\bigset}[2] {\big\{#1\ \big|\ #2\big\}}

\newcommand{\normalpres}[2] {\<#1\ |\ #2\>}
\newcommand{\bigpres}[2] {\big\<#1\ \big|\ #2\big\>}

\newcommand{\normaltuple}[2] {(#1\ |\ #2)}

\newcommand{\lmathbox}[1]{\makebox[0pt][l]{$#1$}}

\newcommand{\cmathbox}[1]{\makebox[0pt]{$#1$}}

\newcommand{\ch}{\,}
\newcommand{\id}{\text{\upshape id}}

\newcommand{\myem}{\bf\boldmath}

\newcommand{\minus}{\smallsetminus}

\newcommand{\daancases}[1] {\left \{\begin{array} {@{}l@{\hspace{4ex}}l@{}} #1 \end{array} \right.}

\psset{dotsize=4pt}

\newcommand{\diagrams}{\psset{arrowsize=2pt 5, 
arrowlength=.6, linewidth=.6pt, nodesep=4pt, arrows=->}}

\newcommand\ra{\rightarrow}
\newcommand\longra{\longrightarrow}
\newcommand\sur{\mathrel{\to\kern-1.8ex\to}}
\newcommand\longsur{\mathrel{\longra\kern-1.8ex\to}}

\newcommand{\equ}{\Leftrightarrow}

\newlength{\daanchop} 

\newcommand{\<}{\langle} 
\renewcommand{\>}{\rangle}

\newcommand{\edoc}{\end{document}}
\newcommand{\be}{\begin{equation}}
\newcommand{\ee}{\end{equation}}
\newcommand{\BOX}{{}\raisebox{-.05ex}{\makebox[1em][r]{$\Box$}}}
\newcommand{\block}{\hfill\BOX\par}
\renewcommand{\qed}{\block}
\newcommand{\col}{\text{\upshape :\ }}

\newcommand{\daandot} {\raisebox{.2ex}{\footnotesize$\circ$}}

\newcommand{\vardot}{\daandot}

\newcommand{\daanlabel} {\makebox[0em][r] {\vardot\hspace{1ex}}}

\newcounter{daana}

\newlength{\daanleftmargin}
\newlength{\daanrightmargin}
\newenvironment{mathlist}{\par%
\begin{list}{\makebox[0em][l]{%
\daanlabel%
\hspace*{-\daanleftmargin}%
\makebox[\columnwidth][r]
  {\refstepcounter{equation}%
  \upshape(\theequation)}}}
{\usecounter{daana}
\setcounter{daana}{\theequation}
\setlength{\rightmargin}{\daanrightmargin}
\setlength{\leftmargin}{\daanleftmargin}
\setlength{\labelwidth}{\daanleftmargin}
\setlength{\labelsep}{0ex}
\setlength{\itemsep}{.7ex}
\setlength{\topsep}{.7ex}
\setlength{\parsep}{0ex}
\setlength{\listparindent}{\parindent}
}}
{\end{list}}

\newenvironment{mathlistn}
{\begin{list}{\daandot}{
\setlength{\itemsep}{0ex}
\setlength{\topsep}{.8ex}
\setlength{\leftmargin}{\daanleftmargin}
\setlength{\labelsep}{.8ex}
\setlength{\labelwidth}{2ex}}}
{\end{list}}

\newcounter{bean}
\newcounter{beann}

\newcommand{\lista}[2]
{\begin{list}{#1}{%
#2%
\setlength{\listparindent}{\parindent}%
\setlength{\leftmargin}{2.3em}%
\setlength{\itemindent}{0mm}%
\setlength{\labelsep}{.4em}%
\setlength{\labelwidth}{10em}%
\setlength{\topsep}{0.7ex}%
\setlength{\itemsep}{0.7ex}%
\setlength{\parsep}{0ex}%
}}

\newenvironment{parts}{\lista{{\upshape(\alph{bean})}}
{\usecounter{bean}}}{\end{list}}

\newcommand{\flushedlista}[2]
{\begin{list}{#1}{%
#2%
\setlength{\listparindent}{\parindent}%
\setlength{\leftmargin}{0mm}%
\setlength{\itemindent}{0mm}%
\setlength{\labelsep}{.4em}%
\setlength{\labelwidth}{-1.7em}%
\setlength{\topsep}{0ex}%
\setlength{\itemsep}{0ex}%
\setlength{\parsep}{0ex}%
#2}}

\newcommand{\daanline}{\hfill\hspace*{0mm}\nopagebreak}

\newcommand{\wild}{}

\newcommand{\thspace}{4mm}

\newtheoremstyle{defi-like}
  {\thspace}
  {\thspace}
  {\rm}
  {}
  {\bfseries}
  {.}
  {.5em}
  {}

\theoremstyle{defi-like}

\newtheorem{defi}[equation]{Definition}

\newtheoremstyle{slanted}
  {\thspace}
  {\thspace}
  {\em}
  {}
  {\bfseries}
  {.}
  {.5em}
  {}

\theoremstyle{slanted}

\newtheorem{prop}[equation]{Proposition}
\newtheorem{conj}[equation]{Conjecture}
\newtheorem{theo}[equation]{Theorem}
\newtheorem{lemm}[equation]{Lemma}
\newtheorem{coro}[equation]{Corollary}
\newtheorem{wildsl}[equation]{\wild}

\newtheoremstyle{roman}
  {\thspace}
  {\thspace}
  {\rm}
  {}
  {\it}
  {.}
  {.5em}
  {}

\theoremstyle{roman}

\newtheorem{wildrm}[equation]{\wild}

\newtheoremstyle{prf}
  {\thspace}
  {\thspace}
  {}
  {}
  {\it}
  {}
  {.5em}
  {}

\theoremstyle{prf}
\newtheorem*{pf}{Proof.\ }

\newtheoremstyle{emptya}
  {\thspace}
  {\thspace}
  {\rm}
  {}
  {\myem}
  {}
  {.5em}
  {}
\theoremstyle{emptya}
\newtheorem*{emptyb}{\wild}

\usepackage{mathpple}\usepackage{charter} 

\newcommand{\daanmath}{\mathbb}
\newcommand{\zz}{\daanmath{Z}}

\newcommand{\rr}{\daanmath{R}}

\usepackage[left=35mm, right=35mm, top=25mm, bottom=25mm, a4paper, includefoot]{geometry}

\begin{document} 
%

\allowdisplaybreaks

\begin{center}
\parbox{.74\textwidth}{%
{\LARGE
\begin{center} 
An asymmetric generalisation of Artin monoids \\
{\normalsize Daan Krammer} \\[-.5ex]
{\normalsize \today} \\[1ex]
\end{center}}}
\end{center}

\begin{abstract}
We propose a slight weakening of the definitions  of Artin monoids and Coxeter monoids. We study one `infinite series' in detail.
\end{abstract}

\newcommand{\arrown}{\mathop{\begin{array}{@{}c@{}} \scriptstyle \smash{\phantom{A}} \\[-1.4ex] \smash{\longra} \\[-1.4ex] \scriptstyle \smash{0} \end{array}}}

\newcommand{\arrownA}{\mathop{\begin{array}{@{}c@{}} \scriptstyle \smash{A} \\[-1.4ex] \smash{\longra} \\[-1.4ex] \scriptstyle \smash{0} \end{array}}}

\newcommand{\arrownM}{\mathop{\begin{array}{@{}c@{}} \scriptstyle \smash{M} \\[-1.4ex] \smash{\longra} \\[-1.4ex] \scriptstyle \smash{0} \end{array}}}

\newcommand{\arrowA}{\stackrel{A}{\longra}}
\newcommand{\arrowwA}{\stackrel{A}{\longsur}}
\newcommand{\arrowM}{\stackrel{M}{\longra}}
\newcommand{\arrowwM}{\stackrel{M}{\longsur}}

\section{Introduction}

This paper begins with a classification of monoids generated by two idempotents such that the ordering of left-division is a lattice ordering.

The result suggests a definition (definition~\ref{vw05}) of a class of monoids which we call AI monoids (A for Artin, I for idempotent). It contains the well-known Artin monoids. 

Every AI monoid comes hand-in-hand with what we call a CI monoid (C for Coxeter, I for idempotent). The twin of an Artin monoid may be called a Coxeter monoid.

An example of an AI monoid is $A_n$ presented by generators $\normalset{p_a}{1\leq a\leq n}$ and relations
\begin{align*}
p_a\, p_b &= p_b\, p_a && \text{if $|a-b|>1$} \\*
p_{a-1}\, p_a\,p_{a-1} &= p_a\,p_{a-1}\,p_a\,p_{a-1} && \text{if $2\leq a\leq n$.}
\end{align*}

The CI monoid $M_n$ of the same type is presented by generators $\normalset{m_a}{1\leq a\leq n}$ and relations
\begin{align*}
m_a\, m_b &= m_b\, m_a && \text{if $|a-b|>1$} \\*
m_{a-1}\, m_a\,m_{a-1} &= m_a\,m_{a-1}\,m_a\,m_{a-1} && \text{if $2\leq a\leq n$} \\
m_{a-1}\, m_a\,m_{a-1} &= m_{a-1}\,m_a\,m_{a-1}\, m_a && \text{if $2\leq a\leq n$} \\*
m_a\, m_a &= m_a && \text{if $1\leq a\leq n$.}
\end{align*}
The monoid $M_n$ appeared earlier in \cite{he}, \cite{o} and \cite{d} as an overarching object in Garside theory; also see section~\ref{vx40}. In \cite{he} and \cite{o} $Q_n$ is the notation for $M_n$.

If a Coxeter group is finite then the corresponding Artin group $A$ is commonly called spherical. Equivalent to this is that any two elements of $A$ have a common right-multiple. Again equivalent is that the corresponding Coxeter monoid $M$ has an element $w_0$, called a sink, such that $x\, w_0\, y=w_0$ for all $x,y\in M$. Again equivalent to this is that the Coxeter {\em monoid\ch} is finite.

We shall show that $M_n$ has a sink. This is proposition~\ref{vx37} and was previously proved in \cite{d} and \cite{he}. On the other hand $M_n$ is infinite if $n\geq 3$ (proposition~\ref{vx41}).

Thus $M_n$ has some properties in common with the spherical Coxeter monoids, some with the nonspherical ones. We feel however that the similarity with the spherical Coxeter monoids is stronger.

As the full class of AI monoids seems beyond reach (even assuming that the corresponding CI monoid has a sink) we decide to focus on the monoids $A_n$ and $M_n$.
Two of our main results, corollaries~\ref{vx17} and \ref{vx42}, are fast solutions to the word problems in $A_n$ and $M_n$. For both monoids we use the shortlex language.

Spherical Artin groups are examples of Garside groups. See \cite{d1} or \cite{d} for Garside theory. Being a Garside group is an elegant and powerful property implying, among others, a fast solution to the word problem. 

Our solution to the word problem for $A_n$ is very different and seems unrelated to Garside properties. Instead we conjecture that $A_n$ is a weak kind of left-Garside monoid, see conjecture~\ref{vx07}. As partial results towards this conjecture we prove that $A_n$ is left-cancellative (proposition~\ref{vx05}) and that it has a Garside element (proposition~\ref{vx06}).

It is known that every Artin monoid $A$ satisfies the so-called cube condition. A closely related property is that if two elements of $A$ have a common right-multiple then they have a least such. AI monoids are not this well-behaved. In section~\ref{vx43} we present an AI monoid which doesn't satisfy the cube condition.

Every Coxeter group comes with a well-known faithful linear representation defined over $\rr$ \cite{h}. In proposition~\ref{vx44} we present a similarly looking linear representation of any CI monoid, with the difference that we make the base ring depend on the Coxeter monoid in question. We don't know if these representations are faithful.

{\bf Acknowledgement.} Many thanks to V.\ Ozornova for pointing out the relevance of the thesis of A.\ Hess \cite{he}.

\section{Monoids generated by two idempotents}

An element $x$ of a monoid is said to be {\em idempotent\ch} if $x^2=x$.

If $a,b$ are elements of a monoid and $n\geq 0$ we write
\[ [a,b;2n]=(ab)^n, \quad [a,b;2n+1]=(ab)^na. \]

A {\em lattice\ch} is an ordered set in which any two elements $x,y$ have a least common upper bound or {\em join\ch} and a greatest common lower bound or {\em meet.}

\begin{prop} \label{vw39} Let $M$ be a monoid generated by two idempotents $a,b$. Let $\leq$ be the relation on $M$ defined by $x\leq y$ if and only if $y\in xM$, in words, $x$ is a left-divisor of $y$. Suppose $1,a,b$ are distinct and neither $a\leq b$ nor $b\leq a$. Then the following are equivalent:
\begin{parts}
\item The relation $\leq$ is an ordering, and a lattice ordering.
\item There are $k,\ell\geq 2$ with $|k-\ell|\leq 1$ satisfying the following. Let $M'$ be the monoid presented by 
\[ M'= \bigpres{A,B}{[A,B;k]=[A,B;k+1]=[B,A;\ell]=[B,A;\ell+1]}. \]
Then there exists an isomorphism $f\col M'\ra M$ such that $f(A)=a$ and $f(B)=b$.
\item After interchanging $a,b$ if necessary there exists $k\geq 2$ such that $M$ admits one of the following presentations:
\begin{align}
& M= \normalpres{a,b}{[a,b;k]=[b,a;k]}, \text{ or} \\
& M= \normalpres{a,b}{[a,b;k]=[a,b;k+1]=[b,a;k+1]}. \label{vw01}
\end{align}
\end{parts}
\end{prop}

The Hasse diagram of $M$ is defined to be the directed graph with vertex set $M$ and which has an arrow labelled $s$ from $x$ to $xs$ whenever $s\in\{a,b\}$ and $x\neq xs$. If (b) holds with $k=3$ and $\ell=4$ then it looks as follows.

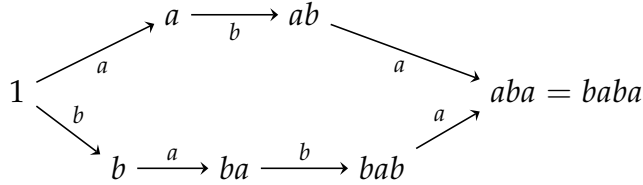
\[
\newcommand{\pscolhookii}{\hspace{-2ex}}
\newcommand{\pscolhookix}{\hspace{-2ex}}
\psmatrix[colsep=2.5ex, rowsep=2.5ex]
&&&& a && ab \\
1 &&&&&&&&&& aba\lmathbox{{}=baba} \\
&&& b && ba && bab 
\endpsmatrix\phantom{{}=baba}
\opt{pics}{%
\diagrams
\everypsbox{\scriptstyle}
\psset{labelsep=2.5pt}
\ncline{2,1}{1,5} \nbput{a}
\ncline{1,5}{1,7} \nbput{b}
\ncline{1,7}{2,11} \nbput{a}
\ncline{2,1}{3,4} \naput{b}
\ncline{3,4}{3,6} \naput{a}
\ncline{3,6}{3,8} \naput{b}
\ncline{3,8}{2,11} \naput{a}
}
\]

\begin{pf} Note that if (b) holds and $\ell=k+1$ then $M$ is presented by (\ref{vw01}). The equivalence (b) $\Leftrightarrow$ (c) is now clear.

Proof of (a) $\Rightarrow$ (b). Since $(M,{\leq})$ is a lattice there exists a join $\Delta$ of $\{a,b\}$. There are $k,\ell\geq 1$ such that
\[ [a,b;k] = \Delta = [b,a;\ell] \]
because $a,b$ are idempotents and $M$ is generated by $a,b$ and $\Delta\in aM$ and $\Delta\in bM$. Choose $k,\ell$ minimal with the above properties. Note $k,\ell\geq 2$ because neither $a\leq b$ nor $b\leq a$.

After interchanging $a,b$ if necessary we may assume $k\leq \ell$.

We have
\begin{align*}
[a,b;k]  \leq [a,b;\ell+1] = a\,[b,a;\ell] = a\Delta= a\,[a,b;k] = [a,b;k]
\end{align*}
so equality holds throughout, proving $[a,b;k]=[a,b;k+1]$. It follows that $\Delta=\Delta a=\Delta b$ (because $a,b$ are idempotents) and therefore $[a,b;\ell]=[a,b;\ell+1]$. 

We shall next prove $\ell \leq k+1$. Suppose to the contrary $\ell\geq k+2$. Put $x=[a,b;\ell-k-2]$ if $k$ is odd and $x=[b,a;\ell-k-2]$ if $k$ is even. Then
\begin{align*}
[b,a;\ell] = b\, [a,b;k+1]\, x = b\, [a,b;k]\, x= b\, [a,b;k] = [b,a;k+1] 
\end{align*}
whence $\ell\leq k+1$ because $\ell$ was chosen minimal. This is a contradiction and proves $\ell \leq k+1$.

We have proved that there exists a unique surjective homomorphism $f\col M'\ra M$ such that $f(A)=a$ and $f(B)=b$. It remains to prove that $f$ is injective. Suppose $x,y\in M$ are distinct with $f(x)=f(y)$. We need to derive a contradiction.

Let $\leq$ be the relation of left division in $M'$ and $0=[A,B;k]=[B,A;\ell]$. Note $0u=u0=0$ for all $u\in M'$.

Suppose first $A\leq x$, $B\leq y$, say, $x=[A,B;p]$ and $y=[B,A;q]$. Then $[a,b;p]=[b,a;q]$. But $a,b$ have a join and $k,\ell$ are minimal so $k\leq p$ and $\ell\leq q$. The definition of $M'$ now implies $x=0=y$, a contradiction.

Suppose next $A\leq x$, $A\leq y$, say, $x=[A,B;p]$, $y=[A,B;q]$. Also assume $p<q$. Then $[a,b;p]=[a,b;q]$. But $\leq$ is an ordering so $[a,b;r]$ is independent of $r$ as long as $p\leq r\leq q$. In particular $[a,b;p]=[a,b;p+1]$. But $a,b$ are idempotents so $f(x)=f(x)a=f(x)b$. So $[a,b;p]=\Delta$. Since $k$ was chosen to be minimal we have $k\leq p<q$. Hence $x=[A,B;p]=[A,B;q]=y$. This is the required contradiction.

Suppose now $1=x$, $A\leq y$. Then $1_M<a=f(A)\leq f(y)=f(x)=1_M$. This   contradicts our assumption that $\leq$ is an ordering.

Up to interchanging $a$ with $b$ or $x$ with $y$ or both this covers all cases. This proves that $f$ is injective and thereby  (a) $\Rightarrow$ (b).

Proof of (a) $\Leftarrow$ (b). Write $\Delta=[a,b;k]=[b,a;\ell]$. Note that $\Delta$ is a sink, that is, $x\Delta y=\Delta$ for all $x,y\in M$. Therefore every element of $M\minus\{1,\Delta\}$ can uniquely be written $[a,b;p]$ ($0<p<k$) or $[b,a;q]$ ($0<q<\ell$).

Conversely, $[a,b;p]\neq\Delta$  and $[b,a;q]\neq\Delta$ if $p<k$ and $q<\ell$ because $|k-\ell|\leq 1$. Therefore the Hasse diagram of $M$ is 
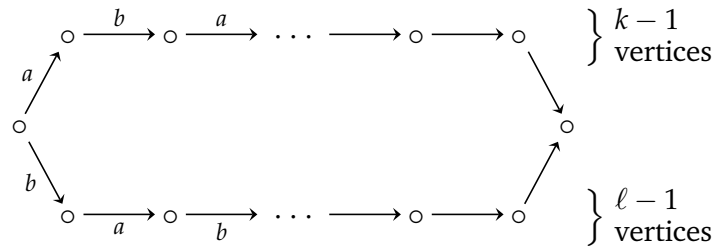
\[
\newcommand{\pscolhookii}{\hspace{-3.5ex}}
\newcommand{\pscolhookvii}{\hspace{-3.5ex}}
\newcommand{\pscolhookviii}{\hspace{-5ex}}
\psmatrix[colsep=5.5ex, rowsep=3ex]
&\circ&\circ&\ \cdots\ \ &\circ&\circ&&\Big\}\ 
\text{\small\begin{tabular}{@{}l@{}}$k-1$\\[-.3ex] vertices\end{tabular}}
\\
\circ &&&&&& \circ & \\
&\circ&\circ&\ \cdots\ \ &\circ&\circ&&\Big\}\ 
\text{\small\begin{tabular}{@{}l@{}}$\ell-1$\\[-.3ex] vertices\end{tabular}}
\opt{pics}{
\diagrams
\everypsbox{\scriptstyle}
\psset{labelsep=2.5pt, nodesep=2.5pt}
\ncline{2,1}{1,2} \naput{a}
\ncline{1,2}{1,3} \naput{b}
\ncline{1,3}{1,4} \naput{a}
\ncline{1,4}{1,5}
\ncline{1,5}{1,6}
\ncline{1,6}{2,7}
\ncline{2,1}{3,2} \nbput{b}
\ncline{3,2}{3,3} \nbput{a}
\ncline{3,3}{3,4} \nbput{b}
\ncline{3,4}{3,5}
\ncline{3,5}{3,6}
\ncline{3,6}{2,7}
}
\endpsmatrix
\]
which proves that $\leq$ is an ordering and a lattice ordering. This finishes the proof of  (a) $\Leftarrow$ (b).\qed\end{pf}


\section{CI monoids and AI monoids}

\begin{defi} \label{vw05} A {\em CI matrix\ch} (C for Coxeter, I for idempotent) consists of a set $S$ and a map $m\col S\times S\ra\zz_{\geq 1}\cup\{\infty\}$ such that:
\begin{mathlistn}
\item $m(a,b)=1$ if and only if $a=b$.
\item $m(a,b)=\infty$ if and only if $m(b,a)=\infty$.
\item $|m(a,b)-m(b,a)|\leq 1$ for all $a,b\in S$.
\end{mathlistn}
With a CI matrix $(S,m)$ we associate the {\em CI monoid\ch} $M$  presented by generating set $S$ and relations
\begin{mathlist}
\item \label{vw02} $a^2=a$ for all $a\in S$.
\item \label{vw03} $[a,b;m(a,b)]=[b,a;m(b,a)]$ whenever $m(a,b)\neq\infty$.
\item \label{vw04} $[a,b;m(a,b)]=[a,b;m(a,b)+1]$ whenever $m(a,b)\neq\infty$.
\end{mathlist}
Moreover we associate an {\em AI monoid\ch} $A$ (A for Artin, I for idempotent) presented by generating set $S$ and relations (\ref{vw03}). 

It is easy to show that the natural map $S\ra M$ is injective. We shall consider $S$ as a subset of $M$ and $A$. Clearly there is a unique homomorphism $A\ra M$ which is the identity on $S$.

A pair $(M',S')$ is called a {\em CI system\ch} if $M'$ is a monoid and there exists an isomorphism $M\ra M'$ (with $M$ as above) taking $S$ to $S'$. Likewise, a pair $(A',S')$ is called an {\em AI system\ch} if $A'$ is a monoid and there exists an isomorphism $A\ra A'$ taking $S$ to~$S'$.

The number $\# S$ is called the {\em rank\ch} of $M$ and $A$.
\qed
\end{defi}

So part (b) of proposition~\ref{vw39} says that $(M,\{a,b\})$ is a CI system of rank~$2$.

Consider definition~\ref{vw05} and suppose that $m$ is symmetric, that is, $m(a,b)=m(b,a)$ for all $a,b\in S$. Then the definition reduces to the following well-known things. Firstly, $(S,m)$ is then known as a Coxeter matrix and $A$ is called an Artin monoid. Likewise we shall call $M$ a Coxeter monoid though this terminology is not common. More commonly studied is the {\em Coxeter group\ch} $W$ which is by definition the quotient of $A$ by the additional relations $a^2=1$ for all $a\in S$ (provided $m$ is symmetric). 

The CI graph or diagram associated with a CI matrix $(S,m)$ is the graph with vertex set $S$ and the following edges:
\begin{mathlistn}
\item If $m(a,b)=m(b,a)>2$ then there is an unoriented edge between $a,b$ labelled $2 m(a,b)=m(a,b)+m(b,a)$.
\item If $m(a,b)+1=m(b,a)$ then there is an arrow from $a$ to $b$ labelled $m(a,b)+m(b,a)$.
\end{mathlistn}
So $m$ is symmetric if and only if all edge labels in the Coxeter graph are even. Warning: If all labels are even then our definition of Coxeter graph differs from the usual one because our labels are twice the usual labels. Labels equal to $6$ are suppressed as usual.

Coxeter groups and Artin monoids have been studied extensively. A good introduction is \cite{h}. Proposition~\ref{vw39} is our main motivation for generalising Artin monoids to AI monoids.

If $M,N$ are monoids then a map $\phi\col M\ra N$ is called an {\em anti-homomor\-phism\ch} if $\phi(xy)=\phi(y)\, \phi(x)$ for all $x,y\in M$.

\begin{lemm} \label{vx23} Let $(M,S)$ be a CI monoid. Assume that no edge label is in $1+4\zz$. That is, $m(a,b)+m(b,a)\not\in 1+4\zz$ for all $a,b\in S$. Then there exists a unique anti-auto\-morphism $\phi$ of $M$ such that $\phi(a)=a$ for all $a\in S$.
\end{lemm}

\begin{pf} We may assume $\# S=2$, say, $S=\{a,b\}$. Write $k=m(a,b)$, $\ell=m(b,a)$. After interchanging $a,b$ if necessary we may also assume $k\leq\ell$.

If $k=\ell$ the result is clear. We are left to consider the case $k<\ell$. Then $\ell=k+1$ and $k$ is odd. By definition $M$ is presented by generating set $\{a,b\}$ and relations $a^2=a$, $b^2=b$ and
\begin{align}
\label{vx25} [a,b;k] &= [b,a;k+1] \\*
\label{vx26} [a,b;k] &= [a,b;k+1] \\
\label{vx27} [b,a;k+1] &= [b,a;k+2].
\end{align}
Note that (\ref{vx27}) is a formal consequence of (\ref{vx25}) and (\ref{vx26}) and can therefore be supressed. The effect of reversing the multiplication is to interchange (\ref{vx25}) and (\ref{vx26}) because $k$ is odd. The result follows.
\qed
\end{pf}

\section{A linear representation for any CI monoid}

The following proposition gives a linear representation of any CI monoid. It looks a bit like the well-known faithful representation of any Coxeter group \cite{h}. We don't know if our representations are faithful.

\begin{prop} \label{vx44} Let $(M,S)$ be a CI system and let $m$ be the associated CI matrix. Let $R$ be the associative ring presented by generators $x_{ab}$ whenever $a,b\in S$ are distinct and relations
\be \label{vx39} [x_{ab},x_{ba};m(a,b)-1]=0 \ee
whenever $a,b\in S$ are distinct. Let $V$ be a free left $R$-module with basis $\normaltuple{e_a}{a\in S}$. Then there exists an $M$-action on $V$ given by
\[ e_a\, a=0, \qquad e_b\, a= e_b+x_{ba}\, e_a \]
whenever $a,b\in S$ are distinct.
\end{prop}

\begin{pf} For $a\in S$ consider the $R$-linear map $T_a\col V\ra V$ defined by
\[ e_a\, T_a=0, \qquad e_b\, T_a= e_b+x_{ba}\, e_a \]
whenever $a,b$ are distinct. Until further notice we shall not use the relations (\ref{vx39}) between the $x_{ab}$. 

We begin by proving that $T_a$ is idempotent. Firstly $e_a\, T_a^2=0=e_a\, T_a$. Moreover if $b\neq a$ then
$ e_b\, T_a^2=(e_b+x_{ba}\, e_a)T_a = e_b\, T_a $
thus proving that $T_a$ is idempotent.

Fix distinct $a,b,c\in S$. For $n\in\zz$ write
\[ \big( a(n), b(n) \big) =  \daancases{(a,b) & \text{if $n$ is even,} \\[.5ex] (b,a) & \text{if $n$ is odd.}} \]

By induction on $n$ we shall prove
\be \label{vx34} e_b\, [T_a,T_b;n] =[x_{ba},x_{ab};n-1]\, e_{a(n)}+[x_{ba},x_{ab};n]\, e_{b(n)} \text{\quad if $n\geq 1$.} \ee
For $n=1$ this is given. If it is true for $n$ then
\begin{align*}
& e_b\, [T_a,T_b;n+1] =  e_b\, [T_a,T_b;n]\, T_{a(n)} \\*& = \big( [x_{ba},x_{ab};n-1]\, e_{a(n)}+[x_{ba},x_{ab};n]\, e_{b(n)} \big) T_{a(n)} \\& = [x_{ba},x_{ab};n]\, e_{b(n)}  T_{a(n)} = [x_{ba},x_{ab};n]\, (e_{b(n)}+x_{b(n),a(n)}\, e_{a(n)}) \\& =  [x_{ba},x_{ab};n]\, e_{b(n)}+[x_{ba},x_{ab};n+1]\, e_{a(n)} \\*& =  [x_{ba},x_{ab};n]\, e_{a(n+1)}+[x_{ba},x_{ab};n+1]\, e_{b(n+1)}.
\end{align*}
This proves (\ref{vx34}).

Since $T_a$, $T_b$ are idempotents
\begin{align} \label{vx35} 
& [T_a,T_b;p+1]-[T_b,T_a;p+1] \\* \notag & \qquad =\big( [T_a,T_b;p]-[T_b,T_a;p]\big)(T_a+T_b-1) \end{align}
for all $p\geq 1$.

By induction on $n$ we shall prove
\begin{align} \label{vx36}
e_c\big([T_a,T_b;n]-[T_b,T_a;n]\big) ={} & x_{ca}[x_{ab},x_{ba};n-1]\, e_{b(n)} \\* \notag {}-{} &  x_{cb}[x_{ba},x_{ab};n-1]\, e_{a(n)} \text{\quad if $n\geq 1$.} \end{align}
It holds for $n=1$ because
\begin{align*}
e_c(T_a-T_b)=(e_c+x_{ca}\, e_a)-(e_c+x_{cb}\, e_b) = x_{ca}\, e_a-x_{cb}\, e_b.
\end{align*}
If it is true for $n-1$ then by (\ref{vx35})
\begin{align*}
& e_c\big([T_a,T_b;n]-[T_b,T_a;n]\big) \\*& =  e_c\big([T_a,T_b;n-1]-[T_b,T_a;n-1]\big) (T_a+T_b-1) \\& = \big( x_{ca}[x_{ab},x_{ba};n-2]\, e_{a(n)} - x_{cb}[x_{ba},x_{ab};n-2]\, e_{b(n)} \big) (T_a+T_b-1) \\& =x_{ca}[x_{ab},x_{ba};n-2]\, e_{a(n)} (T_{b(n)}-1) \\& \,- x_{cb}[x_{ba},x_{ab};n-2]\, e_{b(n)} (T_{a(n)}-1)  \\& =x_{ca}[x_{ab},x_{ba};n-2]\,x_{a(n),b(n)}\,e_{b(n)}  - x_{cb}[x_{ba},x_{ab};n-2]\,x_{b(n),a(n)}\,e_{a(n)}  \\*& =x_{ca}[x_{ab},x_{ba};n-1]\, e_{b(n)}  - x_{cb}[x_{ba},x_{ab};n-1]\, e_{a(n)}.  
\end{align*}
This proves (\ref{vx36}).

We are ready to use the relations (\ref{vx39}) in the ring $R$. Write $k=m(a,b)$ and $\ell=m(b,a)$. 

First suppose $k=\ell$. Write
\[ X=[T_a,T_b;k], \qquad Y=[T_b,T_a;k]. \]
We must prove $X=Y$. Well, (\ref{vx34}) shows that all among $e_aX$, $e_bX$, $e_aY$, $e_bY$ are zero.  Also (\ref{vx36}) shows that $e_c(X-Y)=0$. This settles the case $k=\ell$.

Finally suppose $\ell=k+1$ and write
\[ X=[T_a,T_b;k], \qquad Y=[T_b,T_a;k+1], \qquad Z=[T_a,T_b;k+1]. \]
We must prove $X=Y=Z$. Well, (\ref{vx34}) proves $e_aU=0=e_bU$ for all $U\in\{X,Y,Z\}$. Moreover (\ref{vx36}) shows that $e_c(Y-Z)=0$. 

Applying $[T_a,T_b;k]$ to both sides of the equation $e_c(T_b-1)=x_{cb}\, e_b$ yields
\[ e_c\, \big([T_b,T_a;k+1]-[T_a,T_b;k]\big) = x_{cb}\, e_b\, [T_a,T_b;k] \]
so $e_c(Y-X) = 0$ by (\ref{vx34}). This finishes the case $\ell=k+1$. The proof is complete.\qed
\end{pf}

\section{An AI monoid not satisfying the cube condition} \label{vx43}

It is known that if two elements of an Artin monoid have a common upper bound then they have a join. For AI monoids this is false in general as we shall now show.

Consider the AI monoid $A$ of diagram
\[
\psmatrix[colsep=5ex]
a&b&c.
\endpsmatrix
\opt{pics}{%
\diagrams
\everypsbox{\scriptstyle}
\psset{labelsep=2.5pt}
\ncline{-}{1,1}{1,2} \naput{6}
\ncline{1,2}{1,3} \naput{7}
}
\]
This monoid is presented by
\be A = \bigpres{a,b,c}{aba=bab,\ bcb=cbcb,\ ac=ca}. \label{vx38} \ee
Consider the ordering $\leq$ of left-division on $A$, that is, $x\leq y$ $\Leftrightarrow$ $y=xz$ for some $z$. 

A {\em congruence\ch} on a monoid $N$ is an equivalence relation $\sim$ on $N$ such that there exists a (necessarily unique) structure of monoid on the set $N/{\sim}$ of equivalence classes such that the natural map $N\ra N/{\sim}$ is a homomorphism of monoids.

Let $F$ be the free monoid on $\{a,b,c\}$ and $\sim$ the congruence generated by the relations in (\ref{vx38}), so that $A=F/{\sim}$. For $x\in F$ let $[x]$ denote the equivalence class of $\sim$ containing $x$.

Put $p=[bcb]$, $q=[cabcbab]$ and note
\begin{align*} cabcbab \sim acbcbab \sim abcbab \sim abcaba \sim abacba \sim babcba.
\end{align*}
We have
\[ p=[bcb]=[cbcb], \quad q=[cabcbab]=[babcba] \]
so $p,q$ are two upper bounds of $\{[b],[c]\}$. 

The proof of the following proposition doesn't use any background on Garside theory.

\begin{prop} \daanline
\begin{parts}
\item The set of all words in $a,b,c$ representing $p$ is $\normalset{c^k\,b\,c\,b}{k\geq 0}$.
\item $p$ is a minimal upper bound of $\{[b],[c]\}$. Here minimal means that if $r$ is an upper bound of $\{[b],[c]\}$ with $r\leq p$ then $r=p$.
\item The set of all words representing $q$ is contained in
\begin{align*}
& \bigset{c^k\, b\, c^\ell\, a\, c^m\, b\, c\, b\, a}{k,\ell,m\geq 0} \cup \bigset{c^k\, a\, c^\ell\, b\, a\, c\, b\, a}{k,\ell\geq 0} \\
& \cup \bigset{c^k\, a\, c^\ell\, b\, c\, b\, a\, b}{k,\ell\geq 0} \cup \bigset{c^k\, a\, c^\ell\, b\, c\, a\, b\, a}{k,\ell\geq 0}.
\end{align*}
\item $q$ is not an upper bound of $p$.
\item $\{[b],[c]\}$ has an upper bound but no join.
\end{parts}
\end{prop}

\begin{pf} Parts (a)--(c) are straightforward. By (c) no word for $q$ starts with $bcb$ and so (d) follows. Part (e) follows from (b) and (d).\qed
\end{pf}

There is also a mechanical method for proving that $A$ contains two elements with a common upper bound but without join. To do this one proves that $A$ fails to satisfy the so-called {\em cube condition.} See \cite{d1} or \cite{d} for the necessary background including the $\backslash$ operation. One finds 
\[ \text{$(a\backslash b) \backslash(a \backslash c)=cba$, \qquad $(b\backslash a)\backslash (b\backslash c)=cbab$} \] 
but $cba$, $cbab$ represent distinct elements of $A$.

\section{A CI graph}

From now we shall deal with the CI monoid and the AI monoid of diagram
\be
\begin{array}{c@{}c@{}c@{}c@{}c@{}c@{}c@{}c@{}c}
\circ & \ \stackrel{7}{\longra} \ &
\circ & \  \stackrel{7}{\longra} \  & \ \cdots \ & \  \stackrel{7}{\longra} \  & 
\circ & \  \stackrel{7}{\longra} \  & \circ
\\[-.5ex]
\cmathbox{x_1} &&
\cmathbox{x_2} &&&&&&
\,\,\cmathbox{x_n.}
\end{array} \label{vx24} \ee

Fix a natural number $n$. Let $F_n$ be the free monoid on a set $X_n=\{x_1,\ldots,x_n\}$ of $n$ elements. An element of $F_n$ is called a {\em word} and an element of $X$ a {\em letter.}

\begin{defi} \label{vw89} \daanline \begin{parts}
\item Let $=_B$ be the least congruence on $F_n$ such that
\be 
\label{vw80} 
x_a\, x_b =_B x_b\, x_a \text{\qquad whenever $|a-b|>1$.} \ee
\item Let $=_A$ be the least congruence on $F_n$ containing $=_B$ such that
\begin{align}
x_a\,x_{a-1}\, x_a\,x_{a-1}  =_A x_{a-1}\, x_a\,x_{a-1} \text{\qquad whenever $2\leq a\leq n$.} \label{vw90}
\end{align}
\item Let $=_M$ be the least congruence on $F_n$ containing $=_A$ such that 
\begin{align*}
x_a\, x_a &=_M x_a && \text{for all $a$} \\
x_{a-1}\, x_a\,x_{a-1}\, x_a &= x_{a-1}\, x_a\,x_{a-1} && \text{whenever $2\leq a\leq n$.}
\end{align*}
\end{parts}
An equivalence class with respect to the equivalence relation $=_B$ is called a $B$-class. If $x=_B y$ then we also say that $x$ and $y$ are $B$-equivalent. The $B$-class of $x$ is written $[x]_B$. Likewise for $A$ or $M$ instead of $B$. We write $m_a=[x_a]_M$ and $p_a=[x_a]_A$.

We put 
\[ \text{$A=A_n:=(F_n/{=_A})$, \qquad $M=M_n:=(F_n/{=_M})$.} \] 
\end{defi}

Then $A$ is an AI monoid of diagram (\ref{vx24}) and $M$ a CI monoid of the same diagram.

\section{$M_n$-actions on $X^{n+1}$} \label{vx40}

Let $X$ be a set and write $X^k$ for the Cartesian $k$-th power of $X$. Let $f\col X^2\ra X^2$ be a map. Define maps $f_1,f_2\col X^3\ra X^3$ by $f_1=f\times \id_X$ and $f_2=\id_X\times f$. Assume:
\[ f^2 =f, \qquad f_1\,f_2\,f_1= f_2\,f_1\,f_2\,f_1= f_1\,f_2\,f_1\,f_2. \]
Then there exists an $M_n$-action on $X^{n+1}$ by making $m_a=[x_a]_M$ act as
\[ (\id_X)^{a-1} \times f \times (\id_X)^{n-a}. \]
This simple observation (and the fact that $M_n$ has a sink, see proposition~\ref{vx37}) is at the basis of Garside theory. See \cite{d}, \cite{he}, \cite{o}. This motivates us to focus on $M_n$ and $A_n$.

\section{The diamond lemma}

\begin{lemm}[Diamond lemma] \label{vx08} Let $\ra$ be a relation on a set $S$. Let $\sur$ denote its transitive closure and $\sim$ the equivalence relation generated by $\ra$. Assume:
\begin{mathlist}
\item (Well-founded). There is no infinite sequence $x_1\ra x_2\ra \cdots$ with $x_i\in S$ for all $i$.
\item (Confluence). Let $u,v,w\in S$ and assume $u\ra v$ and $u\ra w$. Then there exists $x\in S$ such that  $v\sur x$ and $w\sur x$. \label{vx09}
\end{mathlist}
An element $v\in S$ is called reduced if there is no $w$ with $v\ra w$. Then every equivalence class for $\sim$ contains a unique reduced element.
\end{lemm}

\begin{pf} See for example \cite[Lemma 1.4.1 and exercise 1.4.2]{c}.\qed
\end{pf}

\section{A rewriting system for $A_n$}

Write $(x_a,x_b]:=x_{a-1}\, x_{a-2}\cdots x_b$ provided $a\geq b$. In particular $(x_a,x_a]=1$. Note also $(x_a,x_b](x_b,x_c]=(x_a,x_c]$.

\begin{defi} Let $\arrownA$ be the least relation on $F_n$ such that
\begin{align}
& x_a\, x_b \arrownA x_b\, x_a \label{vw86}
\end{align}
whenever $a-b\geq 2$ and
\begin{align}
x_{a-1}^{c(1)}\big[ x_{a-2}^{c(2)}\cdots x_{a-b}^{c(b)}\big] (x_a,x_{a-b}] \arrownA \big[  x_{a-2}^{c(2)}\cdots x_{a-b}^{c(b)}\big] (x_a,x_{a-b}] \label{vw87}
\end{align}
whenever $c(i)\geq 1$ for all $i$ and $b\geq 2$. If $u\arrownA v$ then we call $u$ an $A$-standard word. We call the move (\ref{vw86}) a commutation move.
\end{defi}

Note that if $u\arrownA v$ and $u\arrownA w$ then $v=w$. Also, if $u$ and $xuy$ are $A$-standard ($u,x,y\in F_n$) then $x=y=1$.

\begin{defi} \daanline \begin{parts}
\item Let $\arrowA$ be the least relation on $F_n$ containing $\arrownA$ and such that
\[ (u\arrowA v) \Rightarrow (xuy\arrowA xvy) \]
for all $u,v,x,y\in F_n$.
\item We define $\arrowwA$ to be the least transitive relation on $F_n$ containing $\arrowA$.
\end{parts}
\end{defi}

\begin{lemm} \label{vw92} The congruence on $F_n$ generated by $\arrowA$ equals $=_A$.
\end{lemm}

\begin{pf} Let $\sim$ denote the congruence generated by $\arrowA$. 

In (\ref{vw87}) set $b=2$, $c(1)=c(2)=1$. We get
\[ x_{a-1}\, x_{a-2}\,x_{a-1}\, x_{a-2} \sim x_{a-2}\,x_{a-1}\, x_{a-2}. \]
Together with (\ref{vw80}) these generate precisely $=_A$. This proves that $(x =_A y)$ $\Rightarrow$ $(x\sim y)$ for all $x,y\in F_n$. It remains to prove the converse.

Let $x$ be the left-hand side in (\ref{vw87}) and $y$ the right-hand side. Let $P(b)$ denote the statement $x=_Ay$ for all choices of the parameters different from $b$. We will be finished if we can prove $P(b)$ for all $b\geq 2$.

Let $a\in\{2,\ldots,n\}$. By induction on $\ell$ we shall prove
\be \label{vw88} x_{a}\, x_{a-1}^\ell \,x_{a}\, x_{a-1} =_A x_{a-1}^\ell \,x_{a}\, x_{a-1} \ee
for all $\ell\geq 1$. For $\ell=1$ this is (\ref{vw90}). In the following, something in curly brackets is next to be rewritten. Assuming it to be true for $\ell-1$ we find
\begin{align*}
& x_{a}\, \big\{ x_{a-1}^\ell \big\} \,x_{a}\, x_{a-1} \\
&= x_{a}\, x_{a-1}^{\ell-1} \, \big\{ x_{a-1} \,x_{a}\, x_{a-1} \big\} \\
&= \big\{ x_{a}\, x_{a-1}^{\ell-1} \, x_{a}\, x_{a-1} \big\} \,x_{a}\, x_{a-1} && \text{by (\ref{vw90})} \\
&= x_{a-1}^{\ell-1} \, \big\{ x_{a}\, x_{a-1} \,x_{a}\, x_{a-1} \big\} && \text{by the induction hypothesis} \\
&= \big\{ x_{a-1}^{\ell-1} \,  x_{a-1} \big\} \,x_{a}\, x_{a-1} && \text{by (\ref{vw90})} \\
&= x_{a-1}^{\ell}  \,x_{a}\, x_{a-1}.
\end{align*}
We have proved (\ref{vw88}). Using (\ref{vw88}) and an obvious induction on $k$ we find
\[  x_{a}^k\, x_{a-1}^\ell \,x_{a}\, x_{a-1} =_A x_{a-1}^\ell \,x_{a}\, x_{a-1} \]
for all $k\geq 1$ and $\ell\geq 1$. This says that $P(2)$ holds.

We prove $P(b)$ by induction on $b$. Assume $P(b-1)$ and $a-d=b$ and $r(i)\geq 1$ for all $i\in\{a-1,a-2,\ldots,d\}$. 
We simplify notation by writing $e$ instead of $x_e$. We find
\begin{align*}
& (a-1)^{r(a-1)} \cdots d^{r(d)} \big\{ (a,d] \big\} \\*
& = (a-1)^{r(a-1)} \cdots (d+1)^{r(d+1)} \big\{ d^{r(d)}(a,d+2] \big\} (d+2,d] \\
& =_A (a-1)^{r(a-1)} \cdots (d+1)^{r(d+1)} (a,d+2] 
\big\{  d^{r(d)} (d+2,d] \big\} \\
& =_A (a-1)^{r(a-1)} \cdots (d+1)^{r(d+1)} \\*& \qquad\qquad \big\{ (a,d+2] (d+1)\big\} d^{r(d)} (d+2,d] \text{\quad by (\ref{vw88})} \\
& = \big\{ (a-1)^{r(a-1)} \cdots (d+1)^{r(d+1)} (a,d+1] \big\} d^{r(d)} (d+2,d] \\& =_A 
(a-2)^{r(a-2)} \cdots (d+1)^{r(d+1)} \\*& \qquad \big\{(a,d+1]\big\} \,
d^{r(d)} (d+2,d] \text{\quad by the induction hypothesis}  \\& = 
(a-2)^{r(a-2)} \cdots (d+1)^{r(d+1)} (a,d+2] \big\{ (d+1)\, 
d^{r(d)} (d+2,d] \big\} \\& =_A
(a-2)^{r(a-2)} \cdots (d+1)^{r(d+1)} \big\{ (a,d+2]\, d^{r(d)} \big\} (d+2,d] \text{\quad by (\ref{vw88})} \\& =_A
(a-2)^{r(a-2)} \cdots (d+1)^{r(d+1)} \, 
d^{r(d)} \big\{ (a,d+2] (d+2,d] \big\} \\*& =
(a-2)^{r(a-2)} \cdots (d+1)^{r(d+1)} \, 
d^{r(d)} (a,d].
\end{align*}
This proves $P(b-1)\Rightarrow P(b)$ and the proof is complete.
\qed
\end{pf}

\begin{lemm} \label{vw81} The following is the complete list of triples $(q,r,s)$ of nontrivial words such that $qr$ and $rs$ are $A$-standard.
\begin{parts}
\item \label{vw82} A triple $(q,r,s)$ given by
\begin{align*}
q &= x_{a-1}^{r(a-1)} x_{a-2}^{r(a-2)} \cdots x_b^{r(b)} (x_a,x_c] \\[1ex]
r &= (x_c,x_b] \\[1ex]
s &= x_{b-1}^{s(b-1)}  x_{b-2}^{s(b-2)} \cdots x_d^{s(d)} (x_c,x_d]
\end{align*}
whenever 
\[ a\geq c>b\geq d, \qquad a-b\geq 2, \qquad c-d\geq 2 \] 
and $r(i)\geq 1$ for all $i\in\{a-1,a-2,\ldots,b\}$ and $s(j)\geq 1$ for all $j\in\{b-1,b-2,\ldots,d\}$. 
\item \label{vw83} A triple $(q,r,s)$ given by
\begin{align*}
q &= x_c \\[1ex]
r &= x_{a-1} \\[1ex]
s &= \big[x_{a-1}^{r(1)-1}x_{a-2}^{r(2)}\cdots x_{a-b}^{r(b)}\big] (x_a,x_{a-b}] 
\end{align*}
whenever $r(i)\geq 1$ for all $i$ and $b\geq 2$ and $c-a\geq 1$.
\item \label{vw84} A triple $(q,r,s)$ given by
\begin{align*}
q &= \big[x_{a-1}^{r(1)}x_{a-2}^{r(2)}\cdots x_{a-b}^{r(b)}\big] (x_a,x_{a-b+1}] \\[1ex]
r &= x_{a-b} \\[1ex]
s &= x_c
\end{align*}
whenever $r(i)\geq 1$ for all $i$ and $b\geq 2$ and $a-b-c\geq 2$.
\item \label{vw85} A triple $(q,r,s)=(x_a,x_b,x_c)$ where $a-b\geq 2$ and $b-c\geq 2$.
\end{parts}
\end{lemm}

\begin{pf} This is obvious.\qed
\end{pf}

\begin{lemm} \label{vw91}
Let $u,v,w\in F_n$ and assume $u\arrowA v$ and $u\arrowA w$. Then there exists $x\in F_n$ such that  $v\arrowwA x$ and $w\arrowwA x$.
\[
\psmatrix[rowsep=4.5ex, colsep=5.5ex] 
u & v \\
w & x
\endpsmatrix \opt{pics} {\diagrams \everypsbox{\scriptstyle} \psset{nodesep=3pt, labelsep=3pt}
\ncline{1,1}{1,2} \nbput{A}
\ncline{1,1}{2,1} \nbput{A}
\ncline{->>}{1,2}{2,2} \naput{A}
\ncline{->>}{2,1}{2,2} \naput{A}
} 
\]
\end{lemm}

\begin{pf} Throughout the proof we remove the index $A$ from the arrows. 

First suppose there is no overlap, that is, 
\begin{mathlist}
\item \label{vx12} there are $p,q,r,s,t,q',s'\in F_n$ such that $u=pqrst$, $v=pq'rst$, $w=pqrs't$, $q\arrown q'$, $s\arrown s'$. 
\end{mathlist}
Then $x:=pq'rs't$ has the required properties.
\[ \psmatrix[rowsep=4.5ex, colsep=5.5ex] 
pqrst & pq'rst \\
pqrs't & pq'rs't
\endpsmatrix \opt{pics} {\diagrams \psset{nodesep=3pt}
\ncline{1,1}{1,2}
\ncline{1,1}{2,1}
\ncline{->>}{1,2}{2,2}
\ncline{->>}{2,1}{2,2}
} \]
We are left to consider the case of overlap, that is, there are words $p$, $q$, $r$, $s$, $t$, $v_0$, $w_0$ such that
\begin{align*}
& u =p\,q\,r\,s\,t && q\,r \arrown v_0  && r\,s \arrown w_0 \\*
&  r\neq 1   && v = p\,v_0\,s\,t  &&w = p\,q\,w_0\,t.
\end{align*}
We may assume $v\neq w$. It follows that $q\neq 1$ and $s\neq 1$.

We may also assume $p=t=1$.

The possible triples $(q,r,s)$ have been listed in lemma~\ref{vw81}. We shall deal with them one by one.

Suppose first that $(q,r,s)$ is as in lemma~\ref{vw81}(\ref{vw83}). Then, not only can $x_c$ (which is $q$) pass its neighbour $x_{a-1}$ by a commutation move (\ref{vw86}), but it can also go on to pass all remaining letters. This shows
\[ v = v_0\, s\longsur r\, s\, q \longra w_0\, q. \]
Likewise, $q$ can pass all letters in $w_0$ which shows
\[ w = q\, w_0 \longsur w_0\, q. \]
We have shown that $x:=w_0\, q$ has the required properties.

Cases (\ref{vw84}) and (\ref{vw85}) of lemma~\ref{vw81} are similar to case (\ref{vw83}).

It remains to consider case (\ref{vw82}) in lemma~\ref{vw81}. 

Write $e$ instead of $x_e$. In the following, anything between curly brackets is to be rewritten next. On the one hand
\begin{align*}
& \{qr\}s = \big\{ (a-1)^{r(a-1)} \cdots b^{r(b)} (a,b] \big\} (b-1)^{s(b-1)} \cdots d^{s(d)} (c,d] \\*& \ra 
(a-2)^{r(a-2)} \cdots b^{r(b)} (a,b] 
(b-1)^{s(b-1)} \cdots d^{s(d)} (c,d] \\& = 
(a-2)^{r(a-2)} \cdots b^{r(b)} (a,c] \big\{ (c,b] 
(b-1)^{s(b-1)} \cdots d^{s(d)} (c,d] \big\} \\& \ra
(a-2)^{r(a-2)} \cdots b^{r(b)} \big\{ (a,c] (c-1,b] 
(b-1)^{s(b-1)} \cdots d^{s(d)} \big\} (c,d] \\& \sur
(a-2)^{r(a-2)} \cdots b^{r(b)} (c-1,b] 
(b-1)^{s(b-1)} \cdots d^{s(d)} (a,c] (c,d] \\& =
(a-2)^{r(a-2)} \cdots (c-1)^{r(c-1)} \big\{ (c-2)^{r(c-2)} \cdots b^{r(b)} (c-1,b] \big\} \\*& \cdot (b-1)^{s(b-1)} \cdots d^{s(d)} (a,d] =: z.
\end{align*}
Here the last term is abbreviated $z$. On the other hand
\begin{align*}
& q\{rs\} =
(a-1)^{r(a-1)} \cdots b^{r(b)} (a,c] 
\big\{ (c,b] 
(b-1)^{s(b-1)}  \cdots d^{s(d)} (c,d] \big\} \\*& \ra
(a-1)^{r(a-1)} \cdots b^{r(b)}  \big\{ (a,c] 
(c-1,b] (b-1)^{s(b-1)} \cdots d^{s(d)} \big\} (c,d] \\*& \sur
(a-1)^{r(a-1)} \cdots b^{r(b)} 
(c-1,b] (b-1)^{s(b-1)} \cdots d^{s(d)}  (a,d] =: y.
\end{align*}
If $c-b\geq 3$ then
\begin{align*}
& z \ra \big\{ (a-2)^{r(a-2)} \cdots (c-1)^{r(c-1)}  (c-3)^{r(c-3)} \cdots b^{r(b)} \big\} \\*& \cdot (c-1,b] (b-1)^{s(b-1)} \cdots d^{s(d)} (a,d] \sur  (c-3)^{r(c-3)} \cdots b^{r(b)} \\*& \cdot (a-2)^{r(a-2)} \cdots (c-1)^{r(c-1)}   (c-1,b] (b-1)^{s(b-1)} \cdots d^{s(d)} (a,d]\, ; \\[2ex]
& y = 
(a-1)^{r(a-1)} \cdots  (c-1)^{r(c-1)} \big\{ (c-2)^{r(c-2)} \cdots b^{r(b)} (c-1,b] \big\} \\*& \cdot (b-1)^{s(b-1)} \cdots d^{s(d)}  (a,d] \ra \big\{
(a-1)^{r(a-1)} \cdots  (c-1)^{r(c-1)} \\& \cdot   (c-3)^{r(c-3)} \cdots b^{r(b)} \big\} (c-1,b]  (b-1)^{s(b-1)} \cdots d^{s(d)}  (a,d] \\& \sur  (c-3)^{r(c-3)} \cdots b^{r(b)} \\& \big\{
(a-1)^{r(a-1)} \cdots  (c-1)^{r(c-1)}  (c-1,b]  (b-1)^{s(b-1)} \cdots d^{s(d)}  (a,d] \big\} \\& \ra
(c-3)^{r(c-3)} \cdots b^{r(b)} \\*& 
(a-2)^{r(a-2)} \cdots  (c-1)^{r(c-1)}  (c-1,b]  (b-1)^{s(b-1)} \cdots d^{s(d)}  (a,d].
\end{align*}
If $c-b=2$ then
\begin{align*}
& y = \big\{ (a-1)^{r(a-1)} \cdots (b+1)^{r(b+1)} b^{r(b)+1} (b-1)^{s(b-1)} \cdots d^{s(d)}  (a,d] \big\} \\*& \ra (a-2)^{r(a-2)} \cdots (b+1)^{r(b+1)} b^{r(b)+1} (b-1)^{s(b-1)} \cdots d^{s(d)}  (a,d]=z.
\end{align*}
If $c-b=1$ then
\begin{align*}
& y = \big\{ (a-1)^{r(a-1)} \cdots b^{r(b)} (b-1)^{s(b-1)} \cdots d^{s(d)}  (a,d] \big\} \\*& \ra (a-2)^{r(a-2)} \cdots b^{r(b)} (b-1)^{s(b-1)} \cdots d^{s(d)} (a,d]=z.
\end{align*}
This proves the promised result in case (\ref{vw82}) of lemma~\ref{vw81}. The proof is complete.\qed
\end{pf}

\begin{defi} \daanline \begin{parts} 
\item A word $u\in F_n$ is said to be {\em $A$-reduced\ch} if there is no $v$ satisfying $u\arrowA v$.
\item Let $x,y\in F_n$. We say that $x$ is the {\em $A$-reduced form\ch} of $y$ if $x$ is $A$-reduced and $x=_A y$.
\end{parts}
\end{defi}

\begin{theo} \label{vw93} Every $=_A$-class in $F_n$ contains a unique $A$-reduced word.
\end{theo}

\begin{pf} In lemma~\ref{vx08} (the diamond lemma) put $S:=F_n$, $({\ra}):=({\arrowA})$. Then $\sim$ (as defined in the diamond lemma) equals $=_A$ by lemma~\ref{vw92}.

Note that $u\ra v$ implies $\ell(u)>\ell(v)$. Therefore there are no infinite chains $u_1\ra u_2\ra \cdots$. Confluence (\ref{vx09}) is satisfied by lemma~\ref{vw91}. This shows that the assumptions of the diamond lemma are satisfied. The result follows by the diamond lemma.\qed
\end{pf}

\begin{coro} \label{vx17} Consider the AI monoid $A:=(F_n/{=_A})$.\begin{parts}
\item There is a polynomial algorithm computing the $A$-reduced form for a word. 
\item There is a polynomial solution to the word problem in $A$.
\end{parts}
\end{coro}

\begin{pf} (a). Let $u_1\in F_n$ be the input to our algorithm. The algorithm calculates words $u_2,u_3,\ldots,u_n$ such that \[ u_i\arrowA u_{i+1} \] for all $i$, and $u_n$ is $A$-reduced. It is easy to show that each step can be carried out in polynomial time. For all $i$ we have $\ell(u_i)>\ell(u_{i+1})$ so after polynomial time the process terminates, as promised, at some $A$-reduced word $u_n$. The result follows.

(b). This follows immediately from (a) and the fact that every element of $A$ is represented by a unique $A$-reduced word (theorem~\ref{vw93}).\qed
\end{pf}

It would be interesting to know if the methods of this section  apply to the better-known positive braid monoid.

\section{$A_n$ is left-cancellative}

\begin{lemm} \label{vw94} Let $x\in F_n$ be such that $x$ is $A$-standard of length $>2$, that is, $x$ is the left-hand side of {\upshape (\ref{vw87})}. Let $y\in F_n$ be $B$-equivalent to $x$. Then $x$ has the same first letter as $y$, that is, $x=x_a\, u$ and $y=x_a\, v$ for some $a,u,v$.
\end{lemm}

\begin{pf} This is clear.\qed
\end{pf}

Let us call a word {\em $B$-reduced\ch} if it is not of the form $u\, x_a\, x_b\, v$ with $u,v\in F_n$ and $a-b\geq 2$. Clearly, every element of $F_n$ is $B$-equivalent to a unique $B$-reduced word called its {\em $B$-reduced form.}

\begin{prop} \label{vx05} The AI monoid $A_n$ is left-cancellative, that is, if $x,y,z\in F_n$ are such that $xy=_A xz$ then $y=_A z$.
\end{prop}

\begin{pf} Recall that $F_n$ is the free monoid on $X=\{x_1,\ldots, x_n\}$. We may assume $x\in X$.

We may also assume that $y,z$ are $A$-reduced, because otherwise we replace them by their $A$-reduced forms.

Note that $A$-reduced words are $B$-reduced. What does the $B$-reduced form of $xy$ look like? A moment's thought about this question shows that there are words $a,b$ such that $y=ab$ and the $B$-reduced form of $xab$ is $axb$ and such that
\begin{mathlist}
\item $ax=_B xa$. \label{vx29}
\item The letter $x$ doesn't appear in $a$. \label{vx30}
\end{mathlist}
Likewise there are words $c,d$ such that $z=cd$ and the $B$-reduced form of $xcd$ is $cxd$ and such that
\begin{mathlist}
\item $cx=_B xc$. \label{vx31}
\item The letter $x$ doesn't appear in $c$. \label{vx32}
\end{mathlist}
We shall prove that for all $k\geq 1$:
\begin{mathlist}
\item If $axb$ is $A$-reduced then the $A$-reduced form of $x^ky$ is $ax^kb$. \label{vw96}
\item If $axb$ is not $A$-reduced then the $A$-reduced form of $x^ky$ is $ab$. \label{vw95}
\end{mathlist}

Indeed (\ref{vw96}) is immediate. To prove (\ref{vw95}), assume $axb$ is not $A$-reduced. Then there are words $a_1,a_2,b_1,b_2\in F_n$ such that $a=a_1\, a_2$, $b=b_1\, b_2$ and $a_2\, x\, b_1$ is $A$-standard. Using lemma~\ref{vw94} and (\ref{vx29}) and (\ref{vx30}) it follows that $a_2=1$. From (\ref{vw87}) it now follows that
\[ xb_1 \arrowA b_1. \]
Since $\arrowA$ generates $=_A$ as congruence by lemma~\ref{vw92}, we have $xb=_A b$. An obvious induction shows $x^kb=_A b$ and hence $x^ky=_A ax^kb=_A ab$. But $ab$ is $A$-reduced, and we have proved (\ref{vw95}).

Comparison of (\ref{vw96})--(\ref{vw95}) with the analogous statement for $(y,c,d)$ instead of $(x,a,b)$ (in fact the range $k\in\{1,2\}$ is enough) proves that either both $axb$ and $cxd$ are $A$-reduced, or neither is. 

Assume now that $axb$ and $cxd$ are both $A$-reduced. But $axb=_A xy =_A xz =_A cxd$ and an $A$-class doesn't contain more than one $A$-reduced word by theorem~\ref{vw93}. Therefore $axb=cxd$. Also the letter $x$ doesn't appear in $a$ or $c$ by (\ref{vx30}) and (\ref{vx32}). It follows that $a=c$ and $b=d$ and $y=ab=cd=z$. This proves the result if both $axb$ and $cxd$ are $A$-reduced.

Assume finally that $axb$ and $cxd$ are not $A$-reduced. By (\ref{vw95}), $y$ is the $A$-reduced form of $xy$. Likewise $z$ is the $A$-reduced form of $xz$. But $xy=_A xz$ so theorem~\ref{vw93} yields $y=z$. This settles the case where neither $axb$ nor $cxd$ is $A$-reduced. The proof is complete.\qed
\end{pf}

\section{A Garside element in $A_n$}

\begin{defi} \label{vx33} A {\em Garside element\ch} in a monoid $N$ is an element $\Delta\in N$ such that:
\begin{mathlistn}
\item For all $x\in N$ there exist $k\geq 0$ and $y\in N$ such that $xy=\Delta^k$.
\item There exists an endomorphism $\phi$ of $N$ such that $x\Delta=\Delta\phi(x)$ for all $x\in N$.
\end{mathlistn}
\end{defi}

In this section we shall prove that the AI monoid $A_n=(F_n/{=_A})$ has a Garside element.

\begin{defi} We define the elements
\begin{align*} 
Y_n &:= (x_3,x_1]\cdots (x_{n+1},x_1]  \\*[.5ex]
\nabla_n &:= x_1\, Y_n = (x_2,x_1](x_3,x_1]\cdots (x_{n+1},x_1]  
\end{align*}
of $F_n$ and $\Delta_n=[\nabla_n]_A\in A_n$.
\end{defi}

We consider $F_{i-1}$ as a submonoid of $F_i$, for all $i$. Then $\nabla_a\in F_n$ whenever $1\leq a\leq n$. Also $\nabla_n=\nabla_{n-1}(x_{n+1},x_1]$.

\begin{lemm} \label{vw97} We have $x_a\, \nabla_n =_A \nabla_n$ for all $a\in\{2,\ldots,n\}$.
\end{lemm}

\begin{pf} Induction on $n$. For $n=1$ there is nothing to prove. Assume it is true for $n-1$. For $2\leq a\leq n-1$ the induction hypothesis implies
\[ x_a\, \nabla_n=_Ax_a\, \nabla_{n-1}(x_{n+1},x_1] =_A\nabla_{n-1}(x_{n+1},x_1] =_A\nabla_n \] 
thus proving the induction step whenever $2\leq a\leq n-1$. It remains to prove the same for $a=n$. Well,
\begin{align*}
x_n \big\{ \nabla_n \big\} &=_A \big\{ x_n \nabla_{n-2} \big\} (x_{n},x_1] (x_{n+1},x_1] \\*& =_A \nabla_{n-2} \big\{ x_n (x_{n},x_1] (x_{n+1},x_1] \big\} \\*& =_A 
\big\{ \nabla_{n-2} (x_{n},x_1] (x_{n+1},x_1] \big\} && \text{by (\ref{vw87}) and lemma~\ref{vw92}} \\*& =_A \nabla_n.
\end{align*}
This proves the induction step and thereby the lemma.\qed
\end{pf}

\begin{lemm} \label{vw98} We have $x_a\, x_1^r\, \nabla_n =_A  x_1^r\, \nabla_n $ whenever $2\leq a\leq n$ and $0\leq r$.
\end{lemm}

\begin{pf} If $a=2$ then
\begin{align*}
& x_a\, x_1^r \big\{ \nabla_n \big\} \\*& =_A \big\{ x_2\, x_1^r \,x_1 \,x_2 \,x_1 \big\} \,(x_4,x_1]\cdots(x_{n+1},x_1] \\& =_A x_1^r \big\{ x_1 \,x_2 \,x_1 \,(x_4,x_1]\cdots(x_{n+1},x_1] \big\} && \text{by (\ref{vw87}) and lemma~\ref{vw92}} \\*& =_A x_1^r \nabla_n. 
\end{align*}
If $a\geq 3$ then
\begin{align*} 
\big\{ x_a\, x_1^r \big\} \nabla_n &=_A x_1^r \big\{ x_a\,  \nabla_n \big\} \\*& =_A x_1^r\, \nabla_n && \text{by lemma~\ref{vw97}.} \tag*{\BOX}
\end{align*}
\end{pf}

\begin{defi}
Let $\pi\col F_n\ra F_1$ be the homomorphism defined by $\pi(x_1)=x_1$ and $\pi(x_a)=1$ for all $a>1$.
\end{defi}

\begin{lemm} \label{vw99} For all $x\in F_n$ we have $x\, \nabla_n=_A \pi(x)\, \nabla_n$.
\end{lemm}

\begin{pf} For $x\in F_n$, let $k(x):=\ell(x)-\ell(\pi x)$. This is the number of letters in $x$ different from $x_1$. Let $P(n)$ be the statement that the lemma holds whenever $k(x)\leq n$. Then $P(1)$ holds by lemma~\ref{vw98}. 

We prove $P(n)$ by induction on $n$. Assume $P(n-1)$ and let $k(x)=n$. Then we can write $x=yz$ such that $k(y)$ and $k(z)$ are both less than $n$. Then also $k(y\, \pi(z))=k(y)<n$. Using the induction hypothesis we find
\begin{align} \notag x\, \nabla_n & =_A y(z\,\nabla_n) =_A (y\, \pi(z))\nabla_n \\*& =_A \pi(y\,\pi(z))\nabla_n =_A\pi(yz)\nabla_n =_A\pi(x)\nabla_n. \tag*{\BOX} \end{align}
\end{pf}

\begin{lemm} \label{vx02} Let $x\in F_n$. Then there exist $k\geq 0$ and $y\in F_n$ with $xy =_A \nabla_n^k$.
\end{lemm}

\begin{pf} We may assume $n>0$. Then $\pi(\nabla_n)\neq 1$. Therefore there are $z\in F_n$ and $\ell$ such that $\pi(xz)=\pi(\nabla_n^\ell)$. By lemma~\ref{vw99} then
\[ xz\nabla_n =_A \pi(xz)\nabla_n =_A \pi(\nabla_n^\ell)\nabla_n =_A\nabla_n^{\ell+1}. \tag*{\BOX} \]
\end{pf}

\begin{defi} We define an endomorphism $\lambda_n\col F_n\ra F_n$ by $\lambda_n(x_a)=1$ for all $a\in\{2,\ldots,n\}$ and $\lambda_n(x_1)=(x_{n+1},x_n]$.
\end{defi}

\begin{lemm} \label{vx04} For all $x\in F_n$ we have $x\, \nabla_n=_A \nabla_n\, \lambda_n(x)$.
\end{lemm}

\begin{pf} It is clear that we only need to prove this for $\ell(x)=1$. If $x=x_a$ with $a>1$ it follows from lemma~\ref{vw97}. It remains to prove it for $x=x_1$ in which case it states
\[ x_1\, \nabla_n =_A \nabla_n \, (x_{n+1},x_n]. \]
We prove this by induction on $n$. For $n=1$ this is clearly true. Assume it to hold for $n-1$. Then
\begin{align*}
& \nabla_n \, (x_{n+1},x_n] =_A \nabla_{n-1}\, (x_{n+1},x_n]^2 \\*& =_A \nabla_{n-1}\, (x_{n},x_n] (x_{n+1},x_n] && \text{by (\ref{vw87}) and lemma~\ref{vw92}} \\*& =_A x_1\, \nabla_{n-1}\, (x_{n+1},x_n] && \text{by the induction hypothesis} \\*& =_A x_1\, \nabla_{n}.
\end{align*}
This proves the induction step and thereby the lemma.\qed
\end{pf}

\begin{lemm} \label{vx03} If $x,y\in F_n$ are such that $x=_A y$ then $\lambda_n(x)=_A \lambda_n(y)$.
\end{lemm}

\begin{pf} It is enough to prove this if $(x,y)$ is a generator of the congruence $=_A$, that is, $x$ is the left-hand side in (a) or (b) of definition~\ref{vw89} and $y$ the right-hand side. The result is now a simple observation.\qed
\end{pf}

Lemma~\ref{vx03} implies that there exists a unique endomorphism $\phi_n$ of $A_n=(F_n/{=_A})$ such that $\phi_n([x]_A)=[\lambda_n(x)]_A$ for all $x\in F_n$. 

\begin{prop} \label{vx06} The element $\Delta_n$ is a Garside element in $A_n$, with $\phi_n$ playing the role of $\phi$ in definition~{\upshape\ref{vx33}.}
\end{prop}

\begin{pf} This is the content of lemmas~\ref{vx02} and \ref{vx04}.\qed
\end{pf}

We finish with a conjecture.

For $x,y\in A_n$ write $x\leq y$ if and only $y\in xA_n$. This is called the ordering of left division. Note that it is an ordering because $A_n$ is left-cancellative by proposition~\ref{vx05}.

\begin{conj} \daanline \label{vx07}  \begin{parts}
\item The ordered set $(A_n,{\leq})$ is a lattice.
\item Let $p_a$ denote the image of $x_a$ in $A_n$. Let $x\in A_n$. Then $x\leq\Delta_n$ if and only if there exist $z_a\in\<p_2,p_3,\ldots,p_n\>$ for all $a\in\{1,\ldots,n\}$ such that
\[ x= z_n\,(p_1\cdots p_n)\, z_{n-1}\,(p_1\cdots p_{n-1})\cdots z_{2}\, (p_1\,p_2)\, z_1\, p_1. \]
\end{parts}
\end{conj}

A {\em lower semi-lattice\ch} is an ordered set such that any two elements have a meet. 

A weak left-Garside monoid is a monoid with a Garside element and such that the ordering of left-division is a lower semi-lattice. Thus conjecture~\ref{vx07}(a) implies that $A_n$ is a weak left-Garside monoid. The adjective {\em weak\ch} means to remind us that there may be infinitely many left-divisors of $\Delta$, as is the case for $A_n$.

\section{A rewriting system for $M_n$}

Recall the MI monoid $M=M_n=(F_n/{=_M})$ of CI graph (\ref{vx24}). We aim to solve the word problem in this monoid.

\begin{defi} \label{vx19} Let $\arrownM$ be the least relation on $F_n$ such that the following hold.
\begin{parts}
\item $x_a\, x_b\arrownM x_b\, x_a$ whenever $a-b\geq 2$.
\item $(x_a,x_b](x_a,x_b]\arrownM(x_{a-1},x_b](x_a,x_b]$ whenever $a-b\geq 1$. 
In particular, for $a=b$, we have
\[ x_b\, x_b \arrownM x_b. \]
\item \label{vx11} Let $1\leq a\leq b\leq n$. For $i\in\{a+1,\ldots,b\}$ let $y_i$ be an element of the submonoid $\<x_{i+1},x_{i+2},\ldots,x_n\>$ of $F_n$ and $z_i\in \<x_1,\ldots,x_{i-2}\>$. Then we have a rewrite rule
\begin{align}
& x_a(y_{a+1}\, x_{a+1}\, x_a\,z_{a+1} )\cdots (y_{b}\, x_{b}\, x_{b-1}\,z_{b} ) x_{b} \notag \\*
\arrownM {} & x_a(y_{a+1}\, x_{a+1}\, x_a\,z_{a+1} )\cdots (y_{b}\, x_{b}\, x_{b-1}\,z_{b} ). \label{vx18}
\end{align}
In particular, for $a=b$, we have (again)
\[ x_a\, x_a \arrownM x_a. \]
\end{parts}
If $u\arrownM v$ then we call $u$ an $M$-standard word.
\end{defi}

Note that if $u\arrownM v$ and $u\arrownM w$ then $v=w$. Also, if $u$ and $xuy$ are $M$-standard ($u,x,y\in F_n$) then $x=y=1$. As in the case of $A_n$ we define the following.

\begin{defi} \daanline \begin{parts}
\item Let $\arrowM$ be the least relation on $F_n$ containing $\arrownM$ and such that
\[ (u\arrowM v) \Rightarrow (xuy\arrowM xvy) \]
for all $u,v,x,y\in F_n$.
\item We define $\arrowwM$ to be the least transitive relation on $F_n$ containing $\arrowM$.
\end{parts}
\end{defi}

\begin{lemm} \label{vx20} The congruence on $F_n$ generated by $\arrowM$ equals $=_M$.
\end{lemm}

\begin{pf} Let $\sim$ denote the congruence on $F_n$ generated by $\arrowM$. Let $x,y\in F_n$. We must prove
\[ (x \sim y) \equ (x =_M y). \]
The implication $\Leftarrow$ is trivial. In order to prove $\Rightarrow$ we may assume $x\arrownM y$.

If $(x,y)=(x_a\, x_b,x_b\, x_a)$ as in part~(a) of definition~\ref{vx19} then $x=_M y$ is clearly true.

Assume next $x= (x_a,x_b](x_a,x_b]$, $y=(x_{a-1},x_b](x_a,x_b]$ as in part~(b) of definition~\ref{vx19}. Then ${x\arrownA y}$ by (\ref{vw87}) so $x=_A y$ by lemma~\ref{vw92} so $x=_M y$.

Suppose finally that $x$ is the left-hand side in (\ref{vx18}) and $y$ the right-hand side. Note first that the $=_B$-class (commutation class) of $x$ and $y$ doesn't change if we move the $y_i$ all the way to the left and the $z_j$ all the way to the right. Thus we may assume the $y_i$ and $z_j$ to be trivial as we now do.

Let $\rho$ denote the anti-automorphism of $F_n$ defined by $\rho(x_a)=x_a$ for all $a$. By lemma~\ref{vx23} $\rho$ preserves $=_M$. Therefore we need only prove $\rho(x)=_M \rho(y)$. 

Let $u$ be the $B$-reduced form of $\rho(x)$ and $v$ the $B$-reduced form of $\rho(y)$. Then $u=(x_{b},x_a](x_{b},x_a]$ and $v=(x_{b-1},x_a](x_{b},x_a]$. This is precisely a case we've already dealt with. It follows that $u=_M v$ whence $x=_M y$. The proof is finished.\qed
\end{pf}

\begin{lemm} \label{vx10} The following is the complete list of triples $(q,r,s)$ of nontrivial words such that $qr$ and $rs$ are $M$-standard.
\begin{parts}
\item \label{vx14} $(x_a,x_b,x_c)$ whenever $a-b\geq 2$, $b-c\geq 2$.
\item $\big(x_a,x_{b-1},(x_{b-1},x_c](x_b,x_c]\big)$ whenever $a-b\geq 1$, $b-c\geq 1$.
\item $\big((x_a,x_b](x_a,x_{b+1}],x_b,x_c\big)$ whenever $a-b\geq 1$, $b-c\geq 2$.
\item $\big(x_c, x_a, (y_{a+1}\, x_{a+1}\, x_a\,z_{a+1} )\cdots (y_{b}\, x_{b}\, x_{b-1}\,z_{b} ) x_{b} \big)$ when\-ever $c-a\geq 2$ and the notation of {\upshape (\ref{vx11})} holds.
\item \label{vx15} $\big(x_a(y_{a+1}\, x_{a+1}\, x_a\,z_{a+1} )\cdots (y_{b}\, x_{b}\, x_{b-1}\,z_{b} ), x_{b},x_c \big)$\\ whenever $b-c\geq 2$ and the notation of {\upshape (\ref{vx11})} holds.
\item \label{vx13} $\big((x_a,x_b](x_a,x_c],(x_c,x_b],(x_b,x_d](x_c,x_d]\big)$ whenever $a\geq c>b\geq d$.
\item \label{vx16} $\big((x_c,x_a](x_c,x_{a+1}],x_a,(y_{a+1}\, x_{a+1}\, x_a\,z_{a+1} )\cdots (y_{b}\, x_{b}\, x_{b-1}\,z_{b} ) x_{b} \big)$ \\ whenever $c-a\geq 1$ and the notation of {\upshape (\ref{vx11})} holds.
\item \label{vw79} $\big(x_a(y_{a+1}\, x_{a+1}\, x_a\,z_{a+1} )\cdots (y_{b}\, x_{b}\, x_{b-1}\,z_{b} ), x_{b}, (x_{b},x_c](x_{b+1},x_c] \big)$\\ whenever $b-c\geq 0$ and the notation of {\upshape (\ref{vx11})} holds.
\item \label{vw78} Let $1\leq a<b<c\leq n$. For $i\in\{a+1,\ldots,c\}$ let $y_i$ be an element of the submonoid $\<x_{i+1},x_{i+2},\ldots,x_n\>$ of $F_n$ and $z_i$ an element of $\<x_1,\ldots,x_{i-2}\>$. Then we have a triple $(q,r,s)$ with
\begin{align*}
q &= x_a(y_{a+1}\, x_{a+1}\, x_a\,z_{a+1} )\cdots (y_{b}\, x_{b}\, x_{b-1}\,z_{b} )  \\*
r &= x_{b} \\*
s &= (y_{b+1}\, x_{b+1}\, x_{b}\, z_{b+1}) \cdots (y_{c}\, x_{c}\, x_{c-1}\, z_{c}) x_{c}.
\end{align*}
\item \label{vw77} Let $1\leq a<b\leq c\leq n$. For $i\in\{a+1,\ldots,c\}$ let $y_i$ be an element of the submonoid of $\<x_{i+1},x_{i+2},\ldots,x_n\>$ and $z_i$ an element of $\<x_1,\ldots,x_{i-2}\>$. Then we have a triple $(q,r,s)$ with
\begin{align*}
q &= x_a(y_{a+1}\, x_{a+1}\, x_a\,z_{a+1} )\cdots (y_{b-1}\, x_{b-1}\, x_{b-2}\,z_{b-1} )  y_{b}\, x_{b} \\*
r &= x_{b-1}\,  x_{b} \\*
s &= x_{b-1}\, z_{b} (y_{b+1}\, x_{b+1}\, x_{b}\, z_{b+1})  \cdots (y_{c}\, x_{c}\, x_{c-1}\, z_{c}) x_{c}.
\end{align*}
\end{parts}
\end{lemm}

\begin{pf} This is easy.\qed
\end{pf}

\begin{lemm} \label{vx28}
Let $u,v,w\in F_n$ and assume $u\arrowM v$ and $u\arrowM w$. Then there exists $x\in F_n$ such that  $v\arrowwM x$ and $w\arrowwM x$.
\end{lemm}

\begin{pf}
Throughout the proof we remove the index $M$ from the arrows.

If there is overlap (\ref{vx12}) this is proved the same way as in lemma~\ref{vw91}. We are left to consider the case of overlap, that is, there are words $p$, $q$, $r$, $s$, $t$, $v_0$, $w_0$ such that
\begin{align*}
& u =p\,q\,r\,s\,t && q\,r \arrown v_0  && r\,s \arrown w_0 \\*
&  r\neq 1   && v = p\,v_0\,s\,t  &&w = p\,q\,w_0\,t.
\end{align*}
We may assume $v\neq w$. It follows that $q\neq 1$ and $s\neq 1$.

We may also assume $p=t=1$.

The possible triples $(q,r,s)$ have been listed in lemma~\ref{vx10}. We shall deal with them one by one.

\medskip{\myem Cases~(\ref{vx14})--(\ref{vx15}).} In these cases the lemma is readily seen to hold.

\medskip{\myem Case~(\ref{vx13}).} In case~(\ref{vx13}) we write $a$ instead of $x_a$. On the one hand we have
\begin{align*}
& \{qr\}s= \big\{ (a,b](a,c](c,b] \big\} (b,d](c,d] = \big\{ (a,b](a,b] \big\} (b,d](c,d] \\*& \ra (a-1,b](a,b](b,d](c,d] = (a-1,b](a,c] \big\{ (c,d](c,d] \big\} \\& \ra (a-1,b]\big\{(a,c](c-1,d]\big\}(c,d] \sur (a-1,b](c-1,d](a,c](c,d] \\*&
= (a-1,c-1](c-1,b](c-1,b](b,d](a,d] =: y.
\end{align*}
On the other hand
\begin{align*}
& q\{rs\}=(a,b](a,c] \big\{ (c,b](b,d](c,d] \big\} = (a,b](a,c] \big\{ (c,d](c,d] \big\} \\*& \ra (a,b]\big\{(a,c](c-1,d]\big\}(c,d] \sur  (a,b](c-1,d](a,c](c,d] 
\\*& = (a,b](c-1,d](a,d]=:z.
\end{align*}
If $c-b\geq 2$ then
\begin{align*}
& y \ra (a-1,c-1](c-2,b](c-1,b](b,d](a,d] \\*& = \big\{(a-1,c-1](c-2,b]\big\}(c-1,d](a,d] \\*& \sur (c-2,b](a-1,c-1](c-1,d](a,d] = (c-2,b](a-1,d](a,d]\,; \\[2ex]
& z = (a,c-1] \big\{ (c-1,b](c-1,b] \big\} (b,d](a,d] \\& \ra  \big\{(a,c-1](c-2,b]\big\}(c-1,b](b,d](a,d] \\& \sur (c-2,b]\big\{ (a,c-1](c-1,b](b,d] \big\} (a,d] \\*& = (c-2,b] \big\{ (a,d](a,d] \big\} \ra (c-2,b](a-1,d](a,d].
\end{align*}
If $c-b=1$ then
\begin{align*}
z = (a,d](a,d] \ra (a-1,d](a,d] = y.
\end{align*}
Chaining and comparing the above results proves the lemma in case~(\ref{vx13}).

\medskip{\myem Case~(\ref{vx16}).} This case the overlap is untouched, that is, there are words $v_1,w_1\in F_n$ such that $v=v_1\, r\, s$, $w=q\, r\, w_1$. Then $x:=v_1\, r\, w_1$ has the required properties.

\medskip{\myem Case~(\ref{vw79}).} In this case $v=w$ so $x:=v$ has the required properties.

\medskip{\myem Case~(\ref{vw78}).} On the one hand we have
\begin{align*}
& \{qr\}s \ra qs = \big\{ x_a(y_{a+1}\, x_{a+1}\, x_a\,z_{a+1} )\cdots  (y_{c}\, x_{c}\, x_{c-1}\, z_{c})  x_{c} \big\} \\*& \ra x_a(y_{a+1}\, x_{a+1}\, x_a\,z_{a+1} )\cdots  (y_{c}\, x_{c}\, x_{c-1}\, z_{c}).
\end{align*}
On the other hand
\begin{align*}
& q\{rs\} \ra \big\{ x_a(y_{a+1}\, x_{a+1}\, x_a\,z_{a+1} )\cdots (y_b\, x_{b}\, x_{b-1}\,z_b )  \,x_{b} \big\} \\*& (y_{b+1}\, x_{b+1}\, x_{b}\, z_{b+1})\cdots (y_{c}\, x_{c}\, x_{c-1}\, z_{c}) \\*& \ra x_a(y_{a+1}\, x_{a+1}\, x_a\,z_{a+1} )\cdots (y_b\, x_{b}\, x_{b-1}\,z_b ) \\*&  (y_{b+1}\, x_{b+1}\, x_{b}\, z_{b+1})\cdots (y_{c}\, x_{c}\, x_{c-1}\, z_{c}) \\*& = x_a(y_{a+1}\, x_{a+1}\, x_a\,z_{a+1} )\cdots  (y_{c}\, x_{c}\, x_{c-1}\, z_{c}).
\end{align*}
The result follows.

\medskip{\myem Case~(\ref{vw77}).} On the one hand
\begin{align*}
& \{qr\}s \ra q\, x_{b-1}\, s = 
x_a(y_{a+1}\, x_{a+1}\, x_a\,z_{a+1} )\cdots (y_{b-1}\, x_{b-1}\, x_{b-2}\,z_{b-1} ) \\*&  y_{b}\, x_{b}
\big\{ x_{b-1}\, x_{b-1} \big\} z_{b} (y_{b+1}\, x_{b+1}\, x_{b}\, z_{b+1}) \cdots (y_{c}\, x_{c}\, x_{c-1}\, z_{c}) \, x_{c} 
\\& \ra x_a(y_1\, x_{a+1}\, x_a\,z_1 )\cdots  (y_{b-1}\, x_{b-1}\, x_{b-2}\,z_{b-1} ) \,  y_{b}\, x_{b}
\, x_{b-1} \, z_{b} \\& (y_{b+1}\, x_{b+1}\, x_{b}\, z_{b+1}) \cdots  (y_{c}\, x_{c}\, x_{c-1}\, z_{c}) \, x_{c} 
\\& = \big\{ x_a(y_{a+1}\, x_{a+1}\, x_a\,z_{a+1} )\cdots   (y_{c}\, x_{c}\, x_{c-1}\, z_{c}) \, x_{c} \big\}
\\*& \ra x_a(y_{a+1}\, x_{a+1}\, x_a\,z_{a+1} )\cdots (y_{c}\, x_{c}\, x_{c-1}\, z_{c}). 
\end{align*}
On the other hand
\begin{align*}
& q\{rs\} \ra \big\{ x_a(y_{a+1}\, x_{a+1}\, x_a\,z_{a+1} )\cdots (y_{b-1}\, x_{b-1}\, x_{b-2}\,z_{b-1} ) \\*&  y_{b}\, x_{b} \, x_{b-1}\,  x_{b} \big\} \, x_{b-1}\, z_{b} (y_{b+1}\, x_{b+1}\, x_{b}\, z_{b+1})  \cdots (y_{c}\, x_{c}\, x_{c-1}\, z_{c})
\\& \ra  x_a(y_{a+1}\, x_{a+1}\, x_a\,z_{a+1} )\cdots  (y_{b-1}\, x_{b-1}\, x_{b-2}\,z_{b-1} ) \,  y_{b}\, x_{b}  \\& \big\{ x_{b-1}\, x_{b-1} \big\}  z_{b} (y_{b+1}\, x_{b+1}\, x_{b}\, z_{b+1})  \cdots (y_{c}\, x_{c}\, x_{c-1}\, z_{c})
\\& \ra  x_a(y_{a+1}\, x_{a+1}\, x_a\,z_{a+1} )\cdots  (y_{b-1}\, x_{b-1}\, x_{b-2}\,z_{b-1} ) \\&  y_{b}\, x_{b} \, x_{b-1} \, z_{b}  (y_{b+1}\, x_{b+1}\, x_{b}\, z_{b+1})  \cdots (y_{c}\, x_{c}\, x_{c-1}\, z_{c}) \\*& = x_a(y_1\, x_{a+1}\, x_a\,z_1 )\cdots (y_{c}\, x_{c}\, x_{c-1}\, z_{c}).
\end{align*}
This settles case~(\ref{vw77}). The lemma is proved.\qed
\end{pf}

\begin{defi} \daanline \begin{parts} 
\item A word $u\in F_n$ is said to be {\em $M$-reduced\ch} if there is no $v$ satisfying $u\arrowM v$.
\item Let $x,y\in F_n$. We say that $x$ is the {\em $M$-reduced form\ch} of $y$ if $x$ is $M$-reduced and $x=_M y$.
\end{parts}
\end{defi}

\begin{theo} \label{vx21} Every $=_M$-class in $F_n$ contains a unique $M$-reduced word.
\end{theo}

\begin{pf} The proof is the same as for theorem~\ref{vw93}. This time the ingredients are lemmas~\ref{vx20}, \ref{vx28} and \ref{vx08}.\qed
\end{pf}

\begin{coro} \label{vx42} Consider the CI monoid $M:=(F_n/{=_M})$.\begin{parts}
\item There is a polynomial algorithm computing the $M$-reduced form for a word. 
\item There is a polynomial solution to the word problem in $M$.
\end{parts}
\end{coro}

\begin{pf} The proof is the same as for corollary~\ref{vx17}.\qed
\end{pf}

A {\em sink\ch} in a monoid $N$ is an element $0$ such that $x0y=0$ for all $x,y\in N$. 

It is known that a Coxeter monoid is finite if and only if it has a sink. This is false for CI monoids as our next and last two results show.

\begin{prop} \label{vx37} Let $w_0\in M_n$ be the image of $\nabla_n$. Then $x\, w_0\, y=w_0$ for all $x,y\in M_n$.
\end{prop}

\begin{pf} Recall $m_a=[x_a]_M\in M_n$. By lemma~\ref{vw97} we have $m_a\, w_0=w_0$ if $a\geq 2$. Also $m_1^2=m_1$ and $w_0\in m_1\, M_n$ so $m_1\, w_0=w_0$ as well. Since $M_n$ is generated by $\{m_1,\ldots,m_n\}$, we find 
\be \label{vx22} \text{$x\, w_0=w_0$ for all $x\in M_n$.} \ee 

By lemma~\ref{vx23} there exists a unique anti-automorphism $\phi$ of $M_n$ preserving $m_a$ for all $a$. Note that $\phi(w_0)=w_0$. Applying $\phi$ to both sides of (\ref{vx22}) we find $w_0\, y=w_0$ for all $y\in M_n$. The result follows.\qed
\end{pf}

Proposition~\ref{vx37} was earlier proved in \cite[proposition~2.3.14]{he} and \cite{d}.

\begin{prop} \label{vx41} $M_n$ is infinite if $n\geq 3$.
\end{prop}

\begin{pf} Note that $(x_2\, x_1\, x_2\, x_3)^k$ is $M$-reduced for all $k\geq 0$. By theorem~\ref{vx21}(a), they represent distinct elements of $M_n$.\qed
\end{pf}


\end{document}